%% file: cdc2012regularSIAMpaper.tex
\documentclass[letterpaper, 10 pt, conference]{ieeeconf}  

\IEEEoverridecommandlockouts                              
\usepackage{graphics} 
\usepackage{epsfig} 
\usepackage{amsmath} 
\usepackage{amssymb}  
\usepackage{color}
\usepackage{cite}

\def\N{{\mathbb{N}}}
\def\R{{\mathbb{R}}}

\def\Black{\textcolor{black}}
\newcommand{\LG}[1]{\Black{#1}}

\def\state{{x}}
\def\statezero{{\state_{0}}}
\def\statestar{{\state^\star}}
\def\stateset{{X}}
\def\constraintstateset{{\mathbb{\stateset}}}
\def\control{{u}}
\def\controlstar{{\u^\star}}
\def\controlset{{U}}
\def\constraintcontrolset{{\mathbb{\controlset}}}
\def\admissiblecontrolset{{\mathcal{\controlset}}}
\def\feedback{{\mu}}
\def\stateopenloop{{\state_{\control}}}
\def\stagecost{{\ell}}

\def\updateinstant{{\sigma}}
\def\lastupdatedifference{{\varphi}}

\def\stageabbreviation{{\lambda}}

\def\valuefunctionabbreviation{{\nu}}

\def\K{{\mathcal{K}}}

\def\x{{x}}
\def\X{{\stateset}}
\def\Xb{{\constraintstateset}}
\def\u{{\control}}
\def\U{{\controlset}}
\def\Ub{{\constraintcontrolset}}
\def\Uc{{\admissiblecontrolset}}
\def\xu{{\stateopenloop}}
\def\l{{\stagecost}}

\definecolor{karl}{rgb}{.9,.1,.4} 
\definecolor{juergen}{rgb}{.1,.4,.9} 

\newtheorem{theorem}{Theorem}
\newtheorem{definition}[theorem]{Definition}
\newtheorem{proposition}[theorem]{Proposition}

\newtheorem{remark}[theorem]{Remark}
\newtheorem{corollary}[theorem]{Corollary}

\newtheorem{example}[theorem]{Example}

\title{\LARGE \bf
Analysis of Unconstrained Nonlinear MPC Schemes with Time Varying Control Horizon
}

\author{Lars Gr\"{u}ne$^{1}$, J\"{u}rgen Pannek$^{2}$, Martin Seehafer$^{3}$, and Karl Worthmann$^{4}$
\thanks{$^{1}$L.~Gr\"{u}ne is with the Mathematical Institute, University of Bayreuth, 95440 Bayreuth, Germany
        {\tt\small lars.gruene@uni-bayreuth.de}}%
\thanks{$^{2}$J.~Pannek is with the Faculty of Aeronautics, University of the Federal Armed Forces Munich, 85577 Munich/Neubiberg, Germany {\tt\small juergen.pannek@unibw.de}}%
\thanks{$^{3}$M.~Seehafer is with the Munich Reinsurance Company, Divisional Unit: Corporate Underwriting, 80802 Munich, Germany {\tt\small mseehafer@munichre.com}}%
\thanks{$^{4}$K.~Worthmann is with the Mathematical Institute, University of Bayreuth, 95440 Bayreuth, Germany {\tt\small karl.worthmann@uni-bayreuth.de}}%
}

\begin{document}

\maketitle
\thispagestyle{empty}
\pagestyle{empty}

\begin{abstract}
	For nonlinear discrete time systems satisfying a controllability condition, we present a stability condition for model predictive control without stabilizing terminal constraints or costs. The condition is given in terms of an analytical formula which can be employed in order to determine a prediction horizon length for which asymptotic stability or a performance guarantee is ensured. Based on this formula a sensitivity analysis with respect to the prediction and the possibly time varying control horizon is carried out.
\end{abstract}

\section{INTRODUCTION}
By now, model predictive control (MPC) has become a well-established method for optimal control of linear and nonlinear systems, see, e.g., \cite{CamachoBordons2004} and \cite{AllgoewerZheng2000,RawlingsMayne2009}. \LG{The method computes an approximate closed--loop solution to an infinite horizon optimal control problem in the following way: in each sampling interval, based on a measurement of the current state, a finite horizon optimal control problem is solved and the first element (or sometimes also more) of the resulting optimal control sequence is used as input for the next sampling interval(s). This procedure is then repeated iteratively.}

Due to the truncation of the infinite optimization horizon feasibility, stability, and suboptimality issues arise. Suboptimality is naturally discussed with respect to the original infinite horizon optimal control problem, cf. \cite{ShammaXiong1997, NevisticPrimbs1997, GrueneRantzer2008}, but there are different approaches regarding the stability issue. While stability can be guaranteed by introducing terminal point constraints \cite{KeerthiGilbert1988} and \cite{Alamir2006} or Lyapunov type terminal costs and terminal regions \cite{ChenAllgoewer1998,MayneRawlingsRaoScokaert2000}, we focus on a particular stability condition \LG{based on a suboptimality index} introduced in \cite{Gruene2009}, for unconstrained MPC --- that is MPC without modifications such as terminal constraints and costs. \LG{Here, we present a closed formula for this suboptimality index.} This enables us to carry out a detailed sensitivity analysis of this stability criterion with respect to the optimization and the control horizon, i.e., the number of elements of the finite horizon optimal control sequence applied at the plant.

Typically, the length of the optimization horizon predominantly determines the computational effort in each MPC iteration and is therefore considered to be the most important parameter within the MPC method. 
However, suitably choosing the control horizon may lead to enhanced performance estimates and, thus, to significantly shorter optimization horizons. In particular, we prove linear growth of the prediction horizon for appropriately chosen control horizon with respect to a bound on the optimal value function --- an estimate which improves its counterparts given in \cite{GrimmMessinaTeel2005} and \cite{TunaMessinaTeel2006}. Furthermore, we show that MPC is ideally suited in order to deal with networked control systems. 
To this end, the stability proof from \cite{Gruene2009} is extended to time varying control horizons which allows to compensate packet dropouts or non--negligible delays. Here, we show that the corresponding stability condition is not more demanding than its counterpart for so called "classical" MPC for a large class of systems. In addition, \LG{the results in this paper lay the theoretical foundations for MPC algorithms safeguarded by performance estimates obtained for longer control horizons as developed in \cite{PannekWorthmann2011}.}

The paper is organized as follows. In Section \ref{SectionPreliminaries} the problem formulation and the required concepts of multistep feedback laws are given. Then, in Section \ref{SectionResults} a stability condition is derived and analysed with respect to the prediction horizon. In the ensuing Section \ref{SectionTimeVaryingControlHorizon} a stability theorem allowing for time varying control horizon is presented. In order to illustrate our results an example of a nonlinear inverted pendulum on a cart is considered and some conclusions are drawn.

\input{preliminaries}

\input{results}

\input{ResultsTimeVarying}

\input{example}



\section{Conclusions}
We presented a stability condition for MPC without terminal constraints or Lyapunov type terminal costs for nonlinear discrete time systems, which allows to determine a prediction horizon length for which asymptotic stability or a desired guaranteed performance is ensured. Furthermore, we investigated the influence of the prediction and the control horizon on this condition. Suitably choosing the control horizon leads to linear growth of the prediction horizon in terms of the assumed controllability condition. As a consequence, since the prediction horizon predominantly determines the computational costs, computing times can be reduced. In addition, a stability theorem for time varying control horizons was derived. Using symmetry and monotonicity properties, we showed that no additional assumptions were needed in comparison to "classical" MPC.

\addtolength{\textheight}{-12cm}   




\section*{Acknowledgement}

This work was supported by the DFG priority research program 1305 ``Control Theory of Digitally Networked Dynamical Systems'', Grant No. Gr1569/12, and the Leopoldina Fellowship Programme, Grant No. LPDS 2009-36.


\end{document}

%% file: preliminaries.tex
\section{Problem Formulation}\label{SectionPreliminaries}

In this work we consider nonlinear control systems driven by the dynamics
\begin{align}\label{NotationSystemDynamics}
	\state(n+1) = f(\state(n),\control(n))
\end{align}
where $\state$ denotes the state of the system and $\control$ the externally applied control. Both state and control variables are elements of metric spaces $(\X,d_{\X})$ and $(\U,d_{\U})$ which represent the state space and the set of control values, respectively. Hence, our results are also applicable to discrete time dynamics induced by a sampled finite or infinite dimensional system. 
Additionally, state and control are subject to constraints which result in subsets $\constraintstateset \subseteq \stateset$ and $\constraintcontrolset \subseteq \controlset$. Given an initial state $\statezero \in \Xb$ and a control sequence $\control = (u(n))_{n \in I}$, $I = \{0,1,\ldots,N-1\}$ with $N \in \N$ or $I = \N_0$, we denote the corresponding state trajectory by $\stateopenloop(\cdot) = \stateopenloop(\cdot; \statezero)$. Due to the imposed constraints not all control sequences $\control$ lead to admissible solutions. Here, $\admissiblecontrolset^N(x_0)$ denotes the set of all admissible control sequences $\u = (\u(0),\u(1),\ldots,\u(N-1))$ of length $N$ satisfying the conditions $f(\xu(n),u(n)) \in \Xb$ and $\control(n) \in \constraintcontrolset$ for $n = 0, \ldots, N-1$.

We want to stabilize \eqref{NotationSystemDynamics} at a controlled equilibrium $\x^\star$ \LG{and by $\u^\star$ we denote a control value with $f(\x^\star,\u^\star) = \x^\star$.} For given continuous stage costs $\stagecost: \X \times \U \rightarrow \R_0^+$ satisfying $\stagecost(\statestar, \controlstar) = 0$ and $\stagecost(\state, \control) > 0$ for all $\control \in \controlset$ for each $\state \neq \statestar$, our goal is to find a static state feedback law $\stateset \rightarrow \controlset$ which minimizes the infinite horizon cost $J_\infty(\state, \control) = \sum_{n = 0}^\infty \stagecost(\stateopenloop(n), \control(n))$. Since this task is, in general, computationally intractable, we use model predictive control (MPC) instead. Within MPC the cost functional
\begin{align}\label{NotationCostFunctionalFiniteHorizon}
	J_N(\state,\control) := \sum_{n=0}^N \stagecost(\stateopenloop(n;\state),\control(n))
\end{align}
is considered where $N \in \N_{\geq 2}$ denotes the length of the prediction horizon, i.e. the prediction horizon is truncated and, thus, finite. The resulting control sequence itself is also finite. Yet, implementing parts of this sequence, shifting the prediction horizon forward in time, and iterating this procedure ad infinitum yields an implicitly defined control sequence on the infinite horizon. While typically only the first control element of the computed control is applied, cf. \cite{RawlingsMayne2009}, the more general case of multistep feedback laws is considered here. Hence, instead of implementing only the first element at the plant ($m = 1$), $m \in \{1,2,\ldots,N-1\}$ elements of the computed control sequence $u = (u(0),u(1),\ldots,u(N-1))$ are applied. As a result, the system stays in open--loop for $m$ steps. The parameter $m$ is called control horizon.
\begin{definition}[Multistep feedback law]\label{DefinitionMultistepFeedback}
	Let $N \in \N_{\geq 2}$ and $m \in \{1,2,\ldots,N-1\}$ be given. A multistep feedback law is a map $\mu_{N,m}: \X \times \{0,1,\ldots,m-1\} \rightarrow \U$ which is applied according to the rule $\x_{\mu_{N,m}}(0;\x) = \x$,
	\begin{equation*}
		\x_{\mu}(n+1;\x) = f(\x_{\mu}(n;\x),\mu(\x_{\mu}(\varphi(n);\x),n-\varphi(n)))
	\end{equation*}
	with $\mu = \mu_{N,m}$ and $\varphi(n):= \max\{ km | k \in \N_0, km \leq n\}$. 
\end{definition}

For simplicity of exposition, we assume that a minimizer $\control^\star$ of \eqref{NotationCostFunctionalFiniteHorizon} \LG{exists for each $x\in\Xb$ and $N \in \N$. Particularly, this includes the assumption that a feasible solution exists for each $x\in\Xb$. For methods on avoiding this feasibility assumption we refer to \cite{PrimbsNevistic2000} or \cite{GruenePannek2011}.} Using the existence of a minimizer $\controlstar \in \admissiblecontrolset^N(\x)$, we obtain the following equality for the optimal value function defined on a finite horizon
\begin{align}\label{NotationOptimalValueFunction}
	V_N(\state) := \inf_{\control \in \admissiblecontrolset^N(\state)} J_N(\state,\control) = J_N(\x,\u^\star).
\end{align}
Then, the MPC multistep feedback $\mu_{N,m}(\cdot,\cdot)$ is defined by $\mu_{N,m}(\x,n) = \u^\star(n) = \u^\star(n;\x)$ for $n=0,1,\ldots,m-1$. In order to compute a performance or suboptimality index of the MPC feedback $\mu = \mu_{N,m}$, we denote the costs arising from this feedback by
\begin{align*}
	J_\infty^{\feedback}(\state) := \sum_{n=0}^\infty \stagecost \left( \x_\mu(n;\x),\feedback(\x_\mu(\lastupdatedifference(n);\x), n - \lastupdatedifference(n))\right).
\end{align*}

{\em Notation:} throughout this \LG{paper}, we call a continuous function $\rho: \R_{\geq 0} \rightarrow \R_{\geq 0}$ a class $\K_\infty$-function if it satisfies $\rho(0) = 0$, is strictly increasing and unbounded. 

%% file: results.tex
\section{Stability Condition}\label{SectionResults}

In this section we \LG{derive a stability} condition for MPC schemes without stabilizing terminal constraints or costs. To be more precise, we propose a sufficient condition for the relaxed Lyapunov inequality
\begin{equation}\label{TheoremAlphaFormulaRelaxedLyapunovInequality}
	V_N(\x_{\mu}(m;\x)) \leq V_N(\x) - \alpha \sum_{n=0}^{m-1} \l(\x_{\mu}(n;\x),\mu(\x,n)),
\end{equation}
$x \in \Xb$, with $\alpha \in (0,1]$ which, in turn, implies a performance estimate on the MPC closed--loop, cf. \cite{LincolnRantzer2006,GrueneRantzer2008}. We point out that the key assumption needed in this stability condition is always satisfied for a sufficiently large prediction horizon if we suppose that the optimal value function $V_\infty(\cdot)$ is bounded, cf. \cite{AlamirBornard1995, GrimmMessinaTeel2005, JadbabaieHauser2005}. 
In particular, the formula to be deduced allows to easily compute, e.g., a prediction horizon for which stability or a desired performance estimate is guaranteed.
\begin{theorem}\label{TheoremAlphaFormula}
	Let a prediction horizon $N \in \N_{\geq 2}$ and a control horizon $m \in \{1,2,\ldots,N-1\}$ be given. In addition, let a monotone real sequence $\Gamma = (\gamma_0, \gamma_1, \gamma_2, \ldots, \gamma_{N})$, $\gamma_0 = 1$, exist such that the \LG{inequality}
	\begin{equation}\label{TheoremAlphaFormulaControllabilityAssumption}
		V_i(\x) \leq \gamma_i V_1(\x) = \gamma_i \min_{\u \in \Ub: f(\x,\u) \in \Xb} \l(\x,\u) 
		\quad\forall\; \x \in \Xb
	\end{equation}
	holds for all $i \in \{1,2,\ldots,N\}$. Furthermore, \LG{assume that the suboptimality index $\alpha = \alpha_{N,m}$ given by}
	\begin{equation}\label{TheoremAlphaFormulaEq}
		\alpha \hspace*{-0.25mm} = \hspace*{-0.25mm} 1 - \frac {\prod\limits_{i=m+1}^N (\gamma_i - 1) \prod\limits_{i=N-m+1}^N (\gamma_i - 1)}{\left[\prod\limits_{i=m+1}^N \hspace*{-2.5mm} \gamma_i - \hspace*{-2.5mm} \prod\limits_{i=m+1}^N \hspace*{-2.5mm} (\gamma_i - 1) \hspace*{-0.5mm} \right] \hspace*{-1.5mm} \left[\prod\limits_{i=N-m+1}^N \hspace*{-5mm} \gamma_i \hspace*{2mm} - \hspace*{-2mm} \prod\limits_{i=N-m+1}^N \hspace*{-4.5mm} (\gamma_i - 1) \hspace*{-0.5mm} \right]}
	\end{equation}
	\LG{satisfies $\alpha > 0$.} Then, the relaxed Lyapunov Inequality \eqref{TheoremAlphaFormulaRelaxedLyapunovInequality} holds for each $x \in \Xb$ for the feedback law $\mu = \mu_{N,m}$ and the corresponding MPC closed--loop satisfies the performance estimate
	\begin{equation}\label{TheoremAlphaFormulaPerformanceEstimate}
		J_\infty^{\mu_{N,m}}(\x) \leq \frac 1 \alpha V_\infty(\x).
	\end{equation}
	If, in addition, $\mathcal{K}_\infty$-functions $\underline{\eta}$, $\bar{\eta}$ exist such that
	\begin{equation}\label{TheoremAlphaFormulaTechnicalConditions}
		\underline{\eta}(d_{\X}(\x,\x^\star)) \leq V_1(\x) \quad\text{and}\quad V_N(\x) \leq \bar{\eta}(d_{\X}(\x,\x^\star))
	\end{equation}
	hold, then the MPC closed--loop asymptotically converges to $\x^\star$.
\end{theorem}
\proof We sketch the main ideas of the proof and refer for details to \cite{GruenePannekSeehaferWorthmann2010} for the main part and to \cite{Worthmann2011} for the adaptation to our more general setting. 

	Using Bellman's principle of optimality and Condition \eqref{TheoremAlphaFormulaControllabilityAssumption} in order to derive conditions on an open--loop optimal trajectory allows to propose the following optimization problem whose solution yields a guaranteed degree of suboptimality $\alpha$ for the relaxed Lyapunov Inequality \eqref{TheoremAlphaFormulaRelaxedLyapunovInequality}:
	\begin{align*}
		\inf_{\stageabbreviation_0, \ldots, \stageabbreviation_{N-1}, \valuefunctionabbreviation} \frac{\sum_{n=0}^{N-1} \stageabbreviation_n - \valuefunctionabbreviation}{\sum_{n=0}^{m-1} \stageabbreviation_n}
	\end{align*}
	subject to the constraints
	\begin{eqnarray*}
		\sum_{n = k}^{N-1} \stageabbreviation_n & \hspace*{-2.mm} \leq & \hspace*{-2.mm} \gamma_{N-k} \cdot \stageabbreviation_k, \quad k = 0, \ldots, N \hspace*{-0.5mm} - \hspace*{-0.5mm} 2, \\
	\valuefunctionabbreviation - \sum_{n = 0}^{j-1} \stageabbreviation_{n+m} & \hspace*{-2.mm} \leq & \hspace*{-2.mm} \gamma_{N-j} \cdot \stageabbreviation_{j+m}, \quad j = 0, \ldots, N \hspace*{-0.5mm} - \hspace*{-0.5mm} m \hspace*{-0.5mm} - \hspace*{-0.5mm} 1,
	\end{eqnarray*}
	and $\stageabbreviation_0$, \ldots, $\stageabbreviation_{N-1}$, $\valuefunctionabbreviation > 0$. Here, we used the abbreviations $\stageabbreviation_n := \stagecost(\state_{\controlstar}(n), \controlstar(n))$ for a minimizer $\controlstar \in \admissiblecontrolset^N(\state)$ of \eqref{NotationCostFunctionalFiniteHorizon} and $\valuefunctionabbreviation := V_N(\state_{\controlstar}(m))$.

	Within this problem, the constraints represent estimates obtained by using \eqref{TheoremAlphaFormulaControllabilityAssumption} directly or first following an optimal trajectory and, then, making use of \eqref{TheoremAlphaFormulaControllabilityAssumption}. In the next step, this optimization problem is reformulated as a linear program. Then, neglecting some of the imposed inequalities leads to a relaxed linear program whose solution is given by Formula \eqref{TheoremAlphaFormulaEq}. Hence, $\alpha$ from Formula \eqref{TheoremAlphaFormulaEq} is a lower bound for the relaxed Lyapunov Inequality \eqref{TheoremAlphaFormulaRelaxedLyapunovInequality}. 

	If the submultiplicativity condition %
	\begin{equation}\label{TheoremAlphaFormulaSubmultiplicativityCondition}
		\Delta_n \Gamma \cdot \Delta_m \Gamma \geq \Delta_{n+m} \Gamma \quad\text{with}\quad\Delta_i \Gamma := \gamma_{i} - \gamma_{i-1}
	\end{equation}
	is satisfied for all $n,m \in \N$ with $n+m \leq N$ and the given sequence $\Gamma$, Formula \eqref{TheoremAlphaFormulaEq} \LG{actually solves the non-relaxed} problem and, thus, characterizes the desired performance bound even better. \LG{Otherwise, solving the non-relaxed problem may further improve the suboptimality bound $\alpha$.}
\endproof
\begin{remark}
	The main assumption in Theorem \ref{TheoremAlphaFormula} is \LG{Inequality \eqref{TheoremAlphaFormulaControllabilityAssumption}} which is also used in \cite{GrimmMessinaTeel2005,TunaMessinaTeel2006}. However, the performance estimates deduced in these references are more conservative in comparison to the presented technique, cf. \cite{Worthmann2012}. The controllability condition used in \cite{GruenePannekSeehaferWorthmann2010}, i.e. existence of a sequence $(c_n)_{n \in \N_0} \subset \R_{\geq 0}$ such that for each state $x \in \Xb$ an open--loop control $\u_{\x} \in \Uc^\infty(\x)$ exists satisfying
	\begin{equation}\label{NotationControllabilityAssumptionGruene}
		\l(\x_{\u_{\x}}(n;\x),\u_{\x}(n)) \leq c_n V_1(\x)
	\end{equation}
	implies \LG{Inequality \eqref{TheoremAlphaFormulaControllabilityAssumption}} with $\gamma_i := \sum_{n=0}^{i-1} c_n$ but leads, in general, to more conservative estimates, cf. \cite{Worthmann2012}. Note that the methodology proposed in \cite{Gruene2009} allows \LG{to use sequences $\Gamma = (\gamma_i)_{i \in \N}$ depending on the state}. Furthermore, we emphasize that a suitable choice of the stage costs may lead to smaller constants $\gamma_i$, $i \in \{1,2,\ldots,N\}$ and, thus, to improved guaranteed performance, cf. \cite{AltmuellerGrueneWorthmann2010BFG} for an example dealing with a semilinear parabolic \LG{PDE}. 
\end{remark}

We like to mention that Theorem \ref{TheoremAlphaFormula} can be extended to the setting in which an additional weight on the final term is incorporated in the MPC cost functional, i.e.
\begin{equation*}
	J_N(\x_0,\u) := \sum_{n=0}^{N-2} \l(\xu(n),\u(n)) + \omega \l(\xu(N-1),\u(N-1))
\end{equation*}
with $\omega > 1$, cf. \cite[Section 5]{GruenePannekSeehaferWorthmann2010}.

The \LG{availability} of an explicit formula facilitates the analy\-sis of \LG{the performance} estimate $\alpha_{N,m}$ and, thus, allows to draw some conclusions. The first one, \LG{stated formally in the Corollary \ref{cor:convergence}, below,} is that MPC without stabilizing terminal constraints or costs approximates the optimal achievable performance on the infinite horizon arbitrarily well for a sufficiently large prediction horizon $N$ --- independently of the chosen control horizon $m$. For the proof, the concept of an equivalent sequence given in \cite{Worthmann2012} is employed. Then, the argumentation presented in \cite[Corollary 6.1]{GruenePannekSeehaferWorthmann2010} can be used in order to conclude the assertion.
\begin{corollary}
	Let the controllability Condition \eqref{TheoremAlphaFormulaControllabilityAssumption} be satisfied for a monotone bounded sequence $\Gamma = (\gamma_i)_{i \in \N}$. Furthermore, let a control horizon $m \in \N$ be given. Then, the suboptimality estimate $\alpha_{N,m}$, $N \geq \max\{2,m+1\}$, from Formula \eqref{TheoremAlphaFormulaEq} converges to one for $N$ approaching infinity, i.e. $\lim_{N \rightarrow \infty} \alpha_{N,m} = 1$. If, in addition, Condition \eqref{TheoremAlphaFormulaTechnicalConditions} holds, the \LG{MPC closed--loop is asymptotically stable.}
\label{cor:convergence}\end{corollary}

\begin{figure}[thpb]
	\begin{center}
	  \includegraphics[width=6cm]{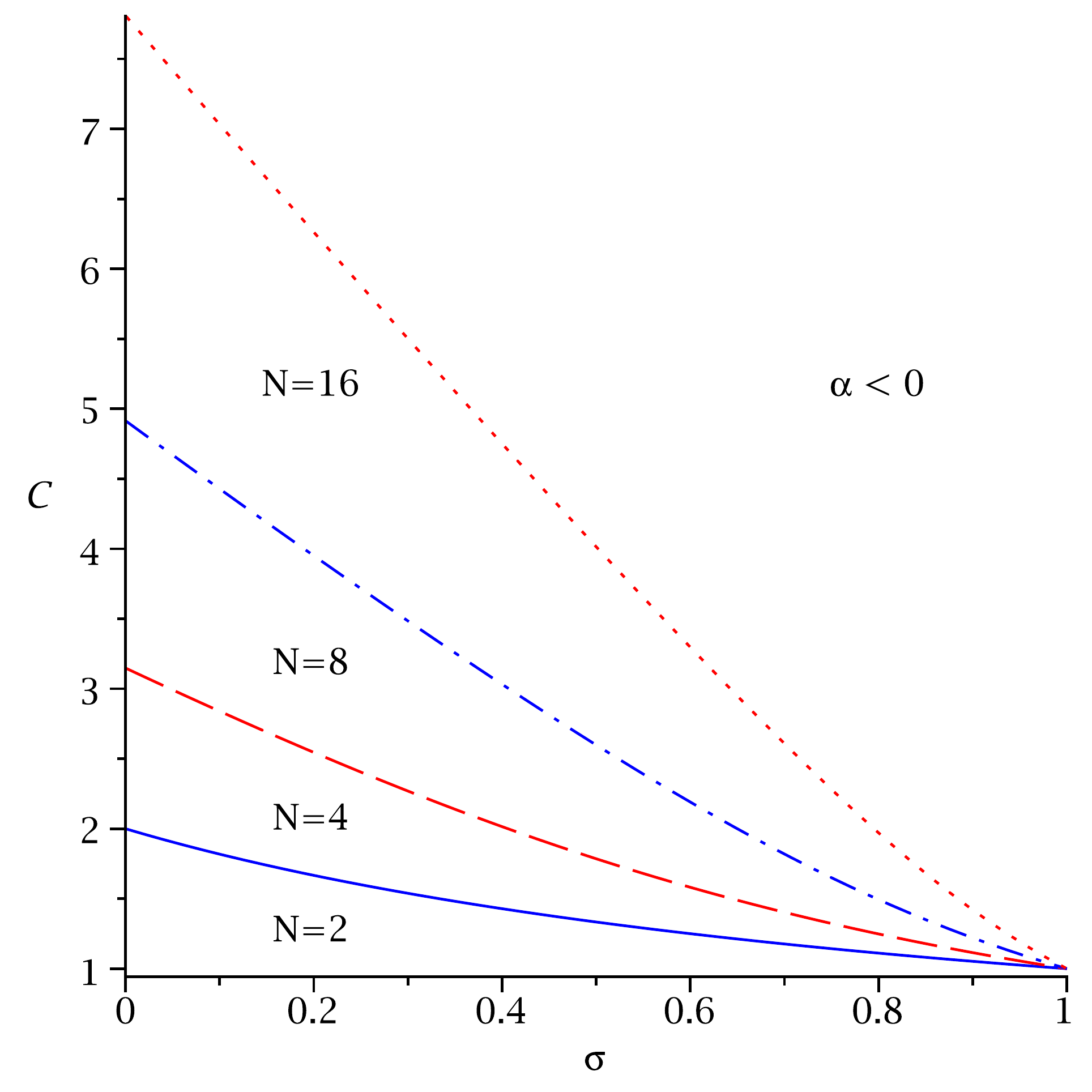}
	  \caption{Parameter pairs $(C,\sigma)$ for function $\beta(r,n) = C \sigma^n r$ such that the corresponding performance bound satisfies the stability condition $\alpha_{N,1} > 0$ depending on the prediction horizon $N$.}
		\label{FigureStabilityRegionDependingOnPredictionHorizon}
	\end{center}
\end{figure}

In order to further elaborate the benefit of \LG{Formula \eqref{TheoremAlphaFormulaEq}}, the following example is considered.
\begin{example}
	Let an exponentially decaying function $\beta(r,n) = C \sigma^n r$ be given. Then, for each prediction horizon $N \in \{2,4,8,16\}$, we determine all parameter combinations $(C, \sigma) \in \R_{\geq 1} \times (0,1)$ such that the stability condition $\alpha_{N,1} \geq 0$ holds with $\gamma_i = C \sum_{n=0}^{i-1} \sigma^n$, cf. Fig. \ref{FigureStabilityRegionDependingOnPredictionHorizon}. Note that this setting corresponds to assuming Condition \eqref{NotationControllabilityAssumptionGruene} with $c_n = C \sigma^n$.
\end{example}

We point out that Fig. \ref{FigureStabilityRegionDependingOnPredictionHorizon} \LG{shows} the different influence of the overshoot $C$ and the decay rate $\sigma$. \LG{Indeed, the figure indicates that for given $N\ge 2$ and $\sigma \in (0,1)$ stability always holds if the overshoot $C > 1$ is sufficiently small. However, for given $N\ge 2$ and overshoot $C>0$ the stability condition may be violated regardless of how $\sigma\in(0,1)$ is chosen. This observation can be proved rigorously using} Formula \eqref{TheoremAlphaFormulaEq}, cf. \cite[Proposition 6.2]{GruenePannekSeehaferWorthmann2010}.

Secondly, Theorem \ref{TheoremAlphaFormula} allows to deduce asymptotic estimates on the minimal prediction horizon length $N$ for which the stability condition $\alpha_{N,m} \geq 0$, $m \in \{1,2,\ldots,N-1\}$, holds --- depending on the sequence $\Gamma = (\gamma_i)_{i \in \N}$ from Condition \eqref{TheoremAlphaFormulaControllabilityAssumption}. Here, one has to keep in mind that the prediction horizon $N$ predominantly determines the required computation time in order to solve the finite horizon optimization problem in each iteration of an MPC algorithm. 

\LG{The next proposition uses a special version of Inequality \eqref{TheoremAlphaFormulaControllabilityAssumption} in which the $\gamma_i$ are independent of $i$. It can be checked, for instance, using an upper bound for the optimal value function $V_\infty$,} cf. \cite[Section 6]{GruenePannekSeehaferWorthmann2010} for a proof.
\begin{proposition}\label{PropositionMinimalStabilizingHorizon}
	Let Condition \eqref{TheoremAlphaFormulaControllabilityAssumption} be satisfied with $\Gamma = (\gamma_i)_{i \in \N}$ with $\gamma_i = M$ for all $i \in \N$.\footnote{Note that the value of $\gamma_1$ is not taken into account in the computation of $\alpha_{N,m}$ from Formula \eqref{TheoremAlphaFormulaEq}. Indeed, \LG{$\gamma_2$ is the first value }contributing to the corresponding suboptimality index.} Then, asymptotic stability of the MPC closed--loop is guaranteed if,
\begin{itemize}
	\item for $m = 1$, the following condition on the optimization horizon is satisfied
		\begin{equation}\label{PropositionMinimalStabilizingHorizon1}
			N \geq 2 + \frac {\ln(M-1)}{\ln(M)-\ln(M-1)}
		\end{equation}
		and, thus, the minimal stabilizing prediction horizon
		\begin{equation}\label{NotationMinimalStabilizingHorizon}
			\hat{N} := \min \{ N : N \in \N_{\geq 2} \text{ and } \alpha_{N,m} \geq 0 \}.
		\end{equation}
		grows asymptotically like $M \ln(M)$ \LG{as $M\to\infty$},
	\item for $m = \lfloor N/2 \rfloor$, one of the following inequalities holds
		\begin{eqnarray}\label{PropositionMinimalStabilizingHorizonN2even}
			\hspace*{-3.5mm} N & \geq & \frac {2\ln(2)}{\ln(M)-\ln(M-1)}, \quad\text{$N$ even}\\
			\hspace*{-3.5mm} N & \geq & \frac{\ln(\frac{2M-1}{M}) \ln(\frac{2M-1}{M-1})}{\ln(M)-\ln(M-1)}, \quad\text{$N$ odd.} \label{PropositionMinimalStabilizingHorizonN2odd}
		\end{eqnarray}
		In this case, the minimal stabilizing Horizon \eqref{NotationMinimalStabilizingHorizon} \LG{grows asymptotically like} $2 \ln(2) M$ \LG{as $M\to\infty$}.
	\end{itemize}
	\LG{By a monotonicity argument, the estimates from this proposition also apply to each sequence $\Gamma = (\gamma_i)_{i \in \N}$ which is bounded by $M$.}
\end{proposition}

\begin{figure}[thpb]
   \includegraphics[width=8cm]{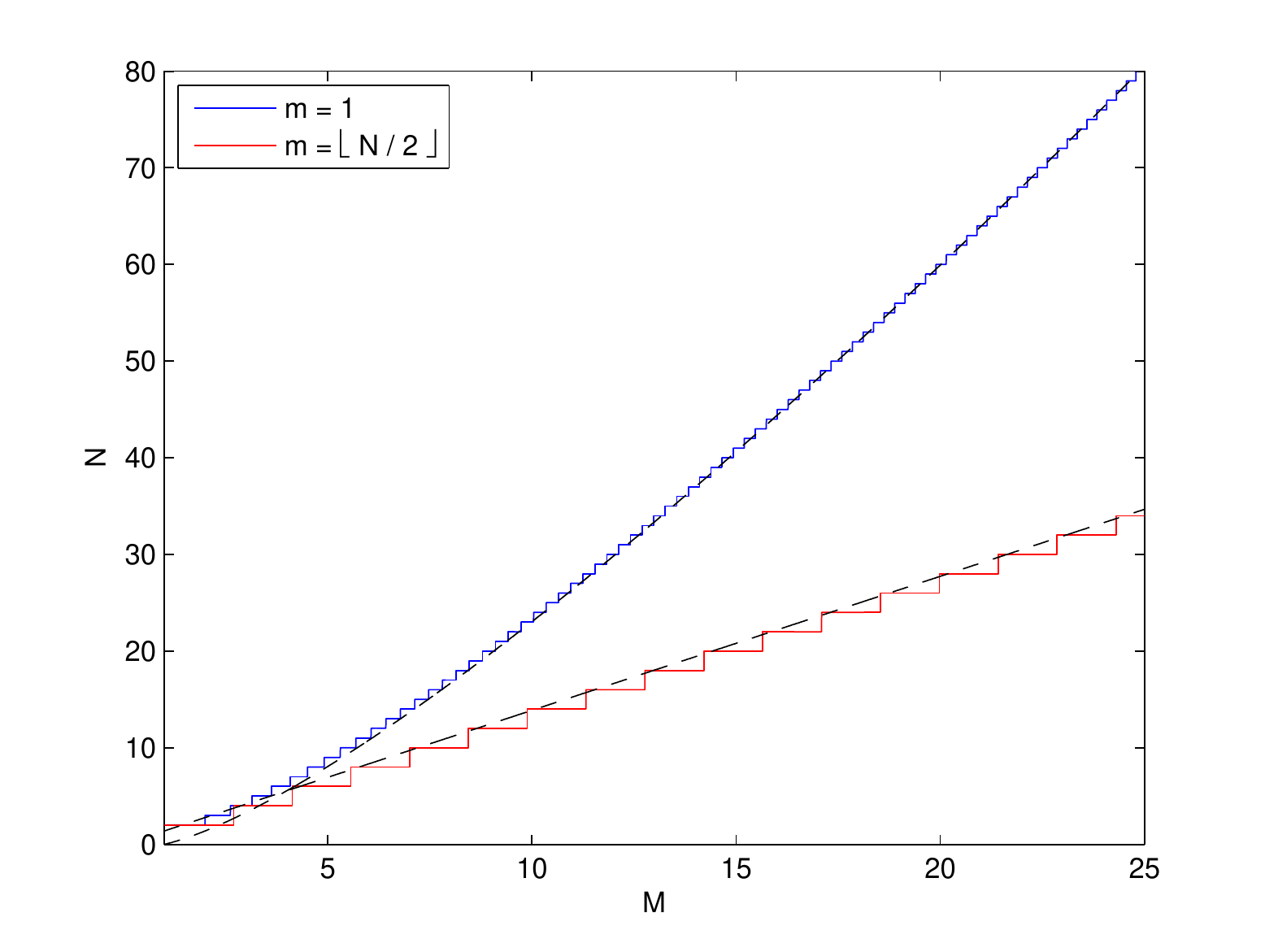}
   \caption{Minimal prediction horizon $N$ for which stability is guaranteed by Theorem \ref{TheoremAlphaFormula} supposing controllability Condition \eqref{TheoremAlphaFormulaControllabilityAssumption} with $(\gamma_i)_{i \in \N}$, $\gamma_i = M$ for all $i \in \N$.}
	\label{FigureMinimalStabilizingHorizon}
\end{figure}

The conclusions of Proposition \ref{PropositionMinimalStabilizingHorizon} are twofold: First, numerical observations from \cite{Gruene2009} are confirmed and the corresponding parameters are precisely determined. Secondly, we emphasize the linear growth of the minimal stabilizing prediction horizon \LG{for $m=\lfloor N/2 \rfloor$.} Hence, the growth for larger control horizons is much slower than for MPC with control horizon $m=1$, cf. Fig. \ref{FigureMinimalStabilizingHorizon}. \LG{This fact will be exploited} \LG{in the following section for both networked control systems and for ``classical'' MPC by designing suitable algorithms.}

%% file: ResultsTimeVarying.tex
\section{Time Varying Control Horizon}\label{SectionTimeVaryingControlHorizon}

In the previous section a stability condition was derived which can also be used in order to ensure a guaranteed performance of the MPC closed--loop. As Proposition \ref{PropositionMinimalStabilizingHorizon} already indicated, employing larger control horizons may improve the corresponding estimates on the required prediction horizon length for which stability can be guaranteed. The following proposition states further properties of \LG{the suboptimality} Bounds \eqref{TheoremAlphaFormulaEq}.
\begin{proposition}\label{PropositionSymmetryMonotonicity}
	Suppose that Condition \eqref{TheoremAlphaFormulaControllabilityAssumption} holds with $\Gamma = (\gamma_i)_{i \in \N}$, $\gamma_i := C \sum_{n=0}^{i-1} \sigma^n$. Here, $C \geq 1$ and $\sigma \in (0,1)$ \LG{denote} overshoot and decay rate of a system which is exponentially controllable in terms of the stage costs. Then, the performance Estimate \eqref{TheoremAlphaFormulaPerformanceEstimate} has the properties:
	\begin{itemize}
		\item symmetry, that is $\alpha_{N,m} = \alpha_{N,N-m}$, and
		\item monotonicity, i.e. $\alpha_{N,m+1} \geq \alpha_{N,m}$ for all $m \in \{1,2,\ldots,\lfloor N/2 \rfloor - 1\}$.
	\end{itemize}
	As a consequence, $\alpha_{N,m} \geq \alpha_{N,1}$ holds for all $m \in \{1,2,\ldots,N-1\}$ \LG{and, in particular, the stability condition $\alpha_{N,m}>0$ holds for arbitrary control horizon $m\ge 2$ if it is satisfied for $m = 1$.}
\end{proposition}
\proof Symmetry follows directly from Formula \eqref{TheoremAlphaFormulaEq}. Contrary to this, showing the claimed monotonicity properties requires a more elaborate technique, cf. \cite[Section 7]{GruenePannekSeehaferWorthmann2010} for a detailed proof.
\endproof

Proposition \ref{PropositionSymmetryMonotonicity} can be exploited in various ways. For instance, in networked control systems the fact $\alpha_{N,m} \geq \alpha_{N,1}$ for all $m \in \{1,2,\ldots,N-1\}$ can be used in order to conclude stability of \LG{a compensation based networked MPC scheme in the presence of packet dropouts or non--negligible delays.} The compensation strategy is \LG{straightforward}: Instead of sending only one control element across the network, an entire sequence is transmitted and buffered at the actuator. If a packet is lost or arrives too late --- that is the packet has not been received by the actuator by the time the first control element of this sequence has to be implemented --- the succeeding element of the current sequence is implemented at the plant which corresponds to incrementing the control horizon $m$. Since it is a priori unknown when and if the next package and, thus, the next sequence of control values arrives at the actuator, the control horizon has to be time varying. Using Theorem \ref{TheoremStabilityTimeVaryingControlHorizon}, \LG{stability can nevertheless} be concluded.

In order to formulate this assertion in a mathematically precise way, the following notation is needed: Let $m^\star \in \{2,\ldots,N-1\}$ be an upper bound for the maximal number of elements of the computed control sequence to be implemented. Then, the transmission times are given by a sequence of control horizons $M = (m_k)_{k \in \N_0}$ with $m^\star \geq m_k \geq 1$. Consequently, in between the $k$th and the $(k+1)$st update of the contol law the system stays in open--loop for $m_k$ steps. Here, we denote the update time instants by $\updateinstant(k) := \sum_{i = 0}^{k-1} m_i$ while $\lastupdatedifference(n) := \max \{ \updateinstant(k) \mid k \in \N_0, \updateinstant(k) \leq n \}$ maps the time instant $n \in \N_0$ to the last update time instant. The corresponding control law is denoted by $\mu_{N,M}$. 
Illustrating these new elements, a control sequence is a sequence
\begin{align*}
	& \feedback(\x_\mu(\updateinstant(k);\x),0),\ldots,\feedback(\x_{\mu}(\updateinstant(k);\x), m_k - 1), \\
	& \quad \feedback(x_{\mu}(\updateinstant(k+1);\x),0), \ldots
\end{align*}
with $\mu = \mu_{N,M}$.
\begin{theorem} \label{TheoremStabilityTimeVaryingControlHorizon}
	Suppose that a multistep feedback law $\feedback_{N, m^\star}: \stateset \times \{ 0, \ldots, m^\star - 1\} \to \controlset$, $m^\star \leq N-1$, and a function $V_N: \stateset \to \R_0^+$ are given. If, for each control horizon $m \in \{1,2,\ldots,m^\star\}$ and each $\state \in \stateset$, we have
	\begin{equation}\label{NotationRelaxedLyapunovInequality}
		V (\x) - V (\x_{\mu}(m;\x)) \geq \alpha \sum_{n = 0}^{m-1} \l(\state_{\mu}(n;\x), \mu(\x,n))
	\end{equation}
	with $\mu = \mu_{N,m^\star}$ for some $\alpha\in(0,1]$, then the estimate $\alpha V_\infty(\x) \leq \alpha V_\infty^{\mu_{N,M}}(\state) \leq V_N(\x)$ holds for all $\x \in \Xb$ and all $M = (m_k)_{k \in \N_0}$ satisfying $m_k \leq m^\star$, $k \in \N_0$. If, in addition, Condition \eqref{TheoremAlphaFormulaTechnicalConditions} is satisfied for $V_N(\cdot)$, asymptotic stability of the MPC closed--loop is ensured.
\end{theorem}

\begin{figure}[thpb]
   \includegraphics[width=8cm]{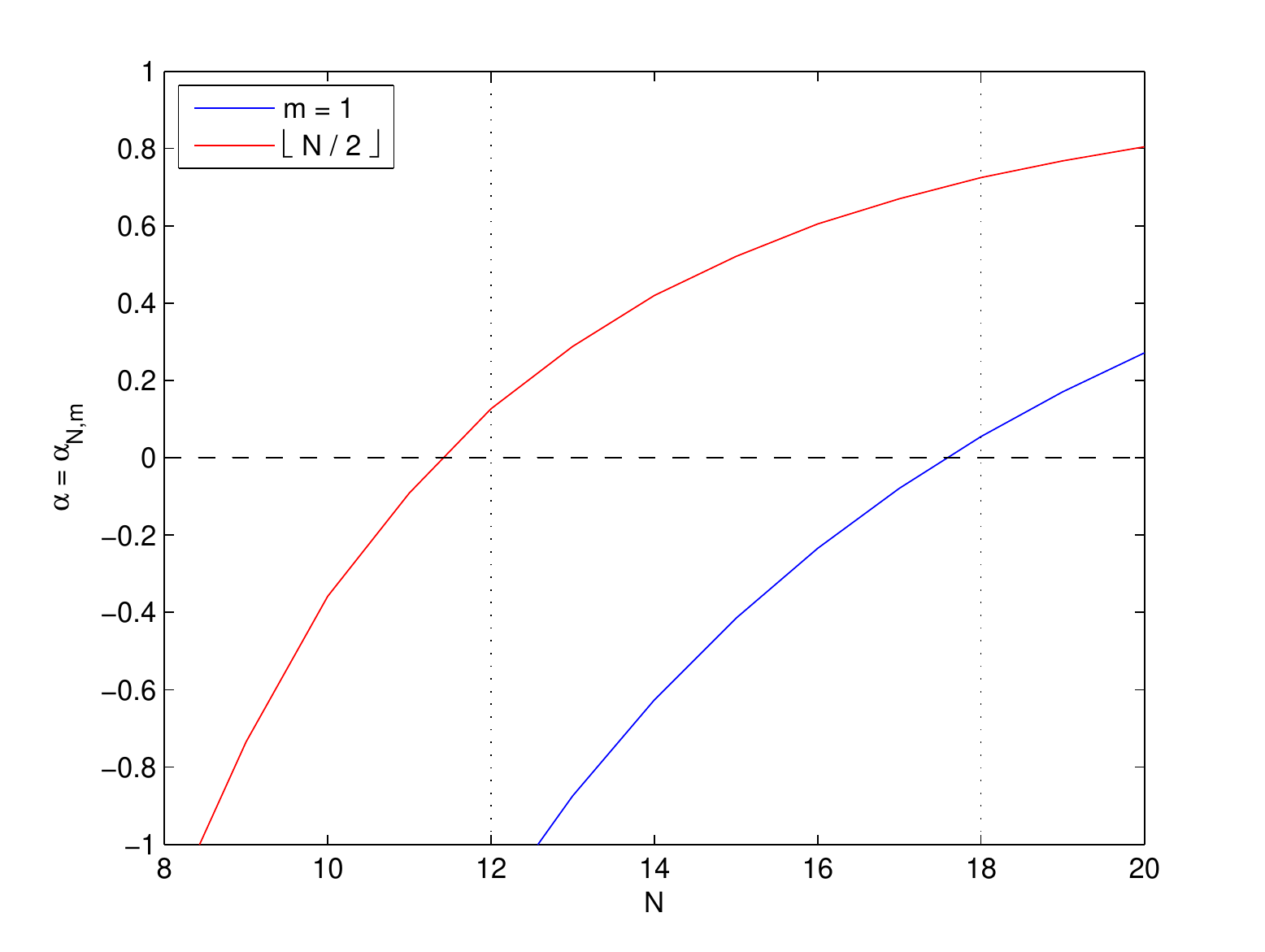}
   \caption{Let Condition \eqref{TheoremAlphaFormulaControllabilityAssumption} be satisfied for the sequence $\Gamma = (\gamma_{i})_{i \in \N}$ with $\gamma_i = C \sum_{n=0}^{i-1} \sigma^n$, $C = 3$, $\sigma = 2/3$. Then, the smallest horizon for which $\alpha_{N,m} \geq 0$ and, thus, stability is guaranteed by Theorem \ref{TheoremAlphaFormula} is $N = 12$. For $m = 1$, even a prediction horizon of length $N = 18$ is required.}\label{FigureNetworkedSystems}
\end{figure}

Theorem \ref{TheoremStabilityTimeVaryingControlHorizon} generalizes its counterpart \cite[Theorem 5.2]{Gruene2009} to time varying control horizon. To this end, the value function $V_N(\cdot)$ was used as a common Lyapunov function, cf. \cite[Theorem 4.2]{GruenePannekSeehaferWorthmann2010}. In order to verify the required assumptions of Theorem \ref{TheoremStabilityTimeVaryingControlHorizon}, our stability condition has to hold for different control horizons $m$, i.e. for each $m \in \{1,2,\ldots,m^\star\}$ which can be checked by Theorem \ref{TheoremAlphaFormula}. However, e.g. for an exponentially controllable system, Proposition \ref{PropositionSymmetryMonotonicity} automatically ensures this condition if it is satisfied for $m = 1$. Hence, the stability condition for time varying control horizons remains the same as for MPC with $m = 1$. Furthermore, we like to point out that increasing the control horizon often enhances the proposed suboptimality bound significantly. In particular, this improvement may lead to a stability guarantee by $\alpha_{N,m} > 0$ although this conclusion cannot be drawn for $m=1$ ($\alpha_{N,1} < 0$), cf. Fig. \ref{FigureNetworkedSystems} and the numerical results shown in Section \ref{SectionExample}.

\LG{Another way to use Proposition \ref{PropositionSymmetryMonotonicity} is described in \cite{PannekWorthmann2011}. There, an algorithm is constructed which employs larger} \LG{control horizons in order to guarantee a desired performance bound. Then, based on an evaluation of the relaxed Lyapunov Inequality \eqref{NotationRelaxedLyapunovInequality}, the MPC loop is closed as often as possible performing a new MPC optimization. This procedure often leads to MPC with $m = 1$, however, safeguarded by the fact that the desired performance can always be ensured --- if necessary --- by enlarging $m$, cf.\ Fig. \ref{FigureNetworkedSystems}. The observed improvement be explained as follows:} Checking the relaxed Lyapunov Inequality \eqref{NotationRelaxedLyapunovInequality} at each time instant is a sufficient but not a necessary condition for \eqref{NotationRelaxedLyapunovInequality} to hold for $m > 1$, i.e. larger control horizons lead to less restrictive conditions.

%% file: example.tex
\section{Example}\label{SectionExample}

\LG{We illustrate our results by computing the $\alpha$-values from the relaxed Lyapunov Inequality \eqref{NotationRelaxedLyapunovInequality} along simulated trajectories in order to compare them with our theoretical findings.} We consider the sampled-data implementation of the nonlinear inverted pendulum on a cart given by the dynamics
\begin{align*}
	\dot{x}_1(t) & = x_2(t)\\
	\dot{x}_2(t) & = -\frac{g}{l}\sin(x_1(t) + \pi) - \frac{k_A}{l} x_2(t) | x_2(t) | \\
	& \qquad - u(t) \cos(x_1(t) + \pi) - k_R \mbox{sgn}(x_2(t))\\
	\dot{x}_3(t) & = x_4(t)\\
	\dot{x}_4(t) & = u(t)
\end{align*}
where $g = 9.81$,  $l = 10$ and $k_R = k_A = 0.01$ denote the gravitation constant, the length of the pendulum and the air as well as the rotational friction terms, respectively. Hence, the discrete time dynamics is defined by $\x(n+1) = \Phi(T;\x(n),\u(n))$. Here, $\Phi(T;\x(n),\u(n))$ represents the solution of the considered differential equation emanating from $x(n)$ with constant control $\u(t) = \u(n)$, $t \in [0,T)$ at time $T$. The goal of our control strategy is to stabilize the upright position $\statestar = (0, 0, 0, 0)$. To this end, we impose the stage cost
\begin{equation*}
	\stagecost(\state(n), \control(n)) := \int_0^T \tilde{\l}(\Phi(t;\x(n),\u(n)),\u(t))\, dt
\end{equation*}
with $\tilde{\l}(\x,\u)$ given by
\begin{align*}
	& 10^{-4} \control^2 + \Big( 3.51 \sin^2 \state_1 + 4.82 \hspace*{0.5mm} \state_2 \sin \state_1  + 2.31 \state_2^2 \\
	+ \hspace*{1.mm} & 0.01 \hspace*{0.5mm} \state_3^2 + 2 \left( (1 - \cos \state_1) (1 + \cos^2 \state_2) \right)^2 + 0.1 \state_4^2 \Big)^2
\end{align*}
with sampling time $T = 0.05$ and prediction horizon $N=70$. Within our computations, we set the tolerance level of the optimization routine and the error tolerance of the differential equation solver to $10^{-6}$ and $10^{-7}$, respectively. Due to the $2 \pi$ periodicity of the stage cost $\stagecost$, we limited the state component $x_1$ to the interval $[-2 \pi + 0.01, 2 \pi - 0.01]$ in order to exclude all equilibria of $\stagecost$ different from $\statestar$. All other state components as well as the control are unconstrained. For our simulations, we used the grid of initial values
\begin{equation*}
      \mathcal{G} := \{\x \in \R^4 | \exists\, i \in \{-1,0,1\}^4: \x = \hat{\x} + 0.05 i \}
\end{equation*}
with $\hat{\x} = (\pi + 1.4, 0, 0, 0)^T$ and computed the suboptimality degree $\alpha_{70,m}$ for constant control horizons $m_k = m$ along the MPC closed loop.\\
Here, we used a startup sequence of $20$ MPC steps with $m = 1$ to compensate for numerical problems within the underlying SQP method. The startup allowed us to compute an initial guess of the optimal open--loop control close to the optimum. During our simulations, we were able to achieve practical stability only, a fact we compensated within our calculations by introducing a truncation region of the stage cost $\stagecost$ using the constant $\varepsilon = 10^{-5}$. The idea of this cut is to take both practical stability regions, that is small areas around the target in which no convergence can be expected, and numerical errors into account, cf. \cite[Theorem 21]{GruenePannek2009} for details. The values of $\alpha_{N,m}$ are computed along the closed--loop trajectory via
\begin{equation}\label{NotationSuboptimalityEstimateAlongTrajectory}
	\alpha_{N,m} = \min_{\x_0 \in \mathcal{G}} \inf_{n \in \{n | \exists k \in \N_0: n = km \}} \alpha_{N,m}(n;\x_0)
\end{equation}
with local degree of suboptimality $\alpha_{N,m}(n)$ given by
\begin{align*}
	\alpha_{N,m}(n;\x_0) = \frac{V_N(\x_{\mu}(n;\x_0)) - V_N(\x_{\mu}(n+m;\x_0))}{\sum\limits_{k=0}^{m-1} (\stagecost(\x_{\mu}(n+k;\x_0),\mu(k;\x_{\mu}(n;\x_0))) - \varepsilon )}
\end{align*}
with $\mu = \mu_{N,m}$ if the denominator of the right hand side is strictly positive and $\alpha_{N,m}(n;\x_0) = 1$ otherwise. Note that $\alpha_{N,m}$ may still become negative if the value function increases along the closed--loop.

\begin{figure}[!ht]
	\begin{center}
		\includegraphics[width=0.48\textwidth]{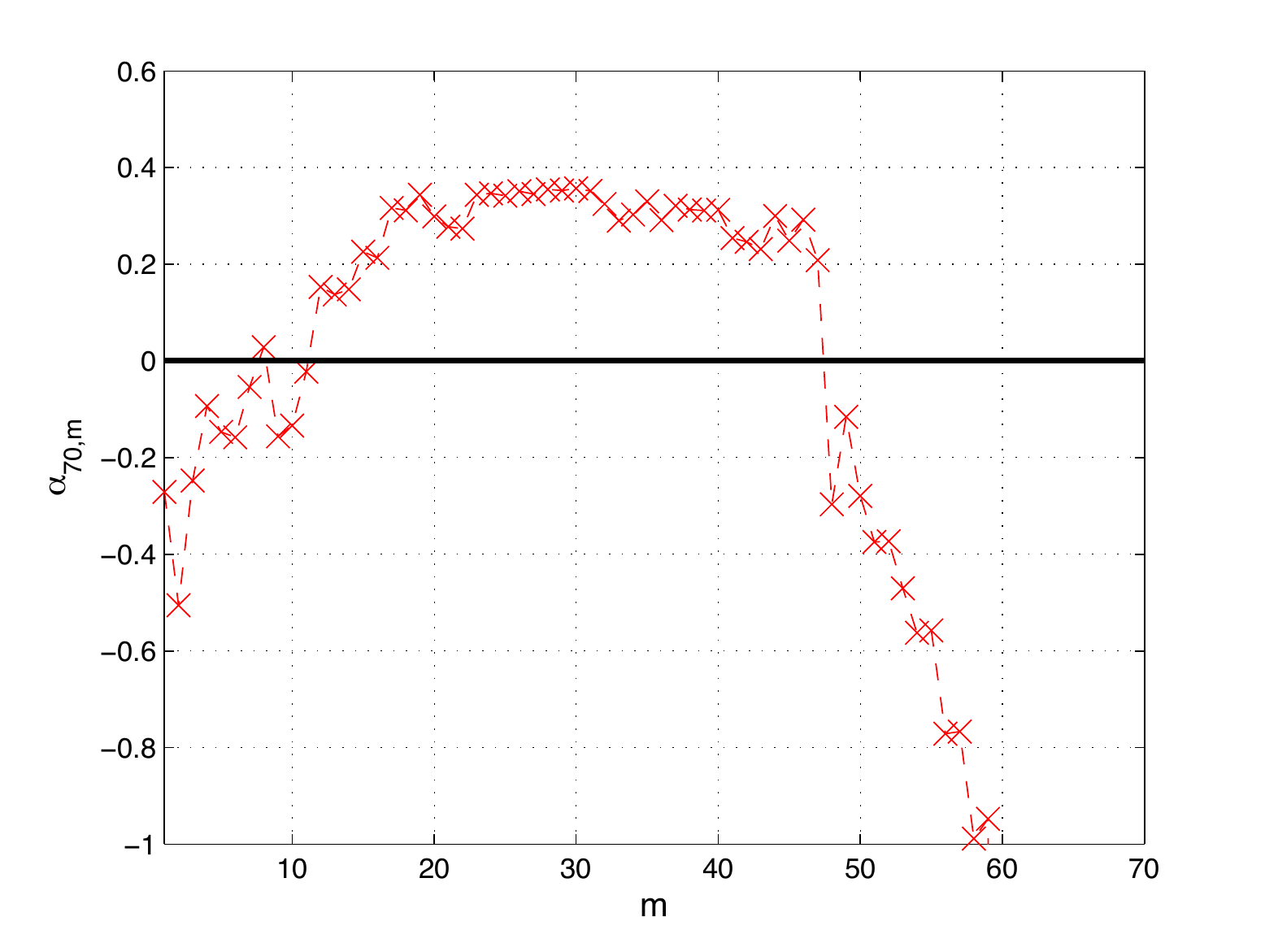}
		\caption{Approximation of $\alpha_{70,m}$, $m \in \{1, \ldots, 70\}$ for the nonlinear inverted pendulum.}
		\label{fig:nonlinear}
	\end{center}
\end{figure}

In Fig.\ \ref{fig:nonlinear}, $\alpha_{N,m}$-values according to \eqref{NotationSuboptimalityEstimateAlongTrajectory} are shown for a variety of control horizons $m$ using the optimization horizon $N = 70$.

While for $m \leq 11$ stability of the closed loop cannot be guaranteed, we obtain $\alpha_{70,m} \geq 0$ for $m \in [12, 47]$. For $m \geq 48$ the values of $\alpha_{70,m}$ are decaying rapidly which may be the result of numerical problems.